\documentclass[graybox]{article}

\usepackage{graphicx}        
\usepackage[bottom]{footmisc}

\usepackage{newtxtext}       %
\usepackage{newtxmath}       
\usepackage{authblk}

\usepackage{url}

\usepackage{geometry}
\usepackage[final]{changes}

\title{A COGENT case study: Supporting Applications with Chombo}

\author[1]{Daniel F. Martin \thanks{DFMartin@lbl.gov} }
\author[2]{Milo Dorr} 
\author[2]{Mikhail Dorf} 
\author[2]{Lee F. Ricketson} 

\affil[1]{Lawrence Berkeley National Laboratory}
\affil[2]{Lawrence Livermore National Laboratory}

\date{}

\begin{document}
\maketitle

\abstract{We present a case study of how a software framework (Chombo) supported the specific needs of a scientific application (COGENT). Since its inception in 2000, the Chombo framework has supported various applications. One example of such support has been the collaboration with the Edge Simulation Laboratory to build the COGENT model.  The specific needs of the COGENT effort required the design and implementation of a set of new capabilities in the Chombo framework, such as higher-order mapped-multiblock discretizations and multi-dimensional code organization.  These capabilities allowed COGENT to develop a unique simulation capability for modeling the edge layers in tokamaks. Once developed, these capabilities were able to support other applications which had similar needs.   }

\section{Introduction}
Practical fusion energy remains a key research goal of the US Department of Energy (DOE) \cite{DOEFusionStrategy} and is one of the NAE's Grand Challenges of Engineering. \cite{NAEreport} The tokamak is a promising design candidate for magnetic confinement fusion energy. \cite{tokamak}  The plasma in a tokamak forms two distinct regions: the core (the inner torus in which the fusion reaction occurs), and the outer edge region (a region of strong gradients, which regulates particle and heat exhaust).  Accurately and efficiently modeling plasmas in the edge regions of tokamaks is essential to understanding the performance the tokamak system, but is extremely challenging due to the complex geometry, strong gradients and wide range of dynamically important length and time scales.  Building such computational models often requires specialized mathematical formulations and their expression in software. In this work, we describe capabilities that the Chombo modeling framework developed in support of the COGENT gyrokinetic modeling code \added{\cite{Dorr2018, Dorf2021}}, and touch on how these capabilities have proved useful to other efforts.

\section{Modeling the Edge Region: The Gyrokinetic Approximation \label{sect:gyrokinetic}}
\begin{figure}
\includegraphics[width=2.5in]{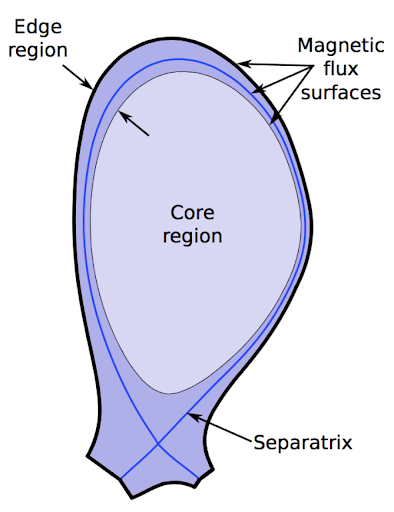}
\caption{Schematic of tokamak edge region -- the separatrix is the outermost magnetic flux surface which is closed. This geometry is also known as an "X-point geometry" due to the crossing of the separatrix -- the specific location is known as the "X-point"}
    \label{fig:tokamakEdgeSchmatic}
\end{figure}

Modeling the edge region of a tokamak presents many numerical difficulties.  The geometry of the edge region itself is complicated, defined on the outside by the device boundaries and centered around a magnetic separatrix, inside of which magnetic field lines are closed, and outside of which, they are open (Figure \ref{fig:tokamakEdgeSchmatic}).
Transport in this region is highly anisotropic, with characteristic lengths along magnetic field lines being about three orders of magnitude larger than that in the perpendicular direction. Due to weak collisionality just inside the magnetic separatrix, a kinetic description of a tokamak edge plasma is needed. Moreover, the presence of steep edge gradients generates large deviations from local thermodynamic equilibrium.  This necessitates a so-called full-F simulation model \added{(such models use both continuum \cite{Dorf2021} or Particle-In-Cell (PIC) \cite{Chang2017} discretizations)}, in contrast to the simplified and more efficient delta-F models \added{(\cite{Candy2003} and \cite{Chen2007})} available in the core due to the proximity to local thermodynamic equilibrium.

Strong plasma magnetization enables the use of the gyrokinetic approximation, in which averaging over fast particle gyromotion is performed and a full 6D kinetic particle distribution can be accurately represented by a 5D gyro-averaged distribution function. Still, a 5D strongly anisotropic transport problem has to be solved in a X-point geometry. 

Although particle-in-cell codes are less sensitive to geometrical features and have been used for edge plasma modeling \cite{Chang2017}, the signal-to-noise ratio decreases only as $1/\sqrt{N}$ in full-F simulations (in contrast to $1/N$ dependence for the delta-f calculations suitable for the core region), which motivates development of a continuum approach.

\added{In a full-F approach, we simultaneously model large-scale/large-amplitude background dynamics and small-scale/small-amplitude turbulence. High numerical accuracy (in the form of higher-order discretizations) is particularly important to ensure that numerical (truncation) errors from the large-amplitude background do not overwhelm small-amplitude turbulence processes. }

\added{Additionally, c}\deleted{C}omputationally efficient continuum modeling of strongly anisotropic transport requires the use of meshing which is aligned with magnetic field lines and magnetic flux surfaces. However, such a meshing approach poses significant challenges due to diverging metric factors at the X-point. Developed by the Edge Simulation Laboratory collaboration, the finite-volume \deleted{gyrokientic}\added{gyrokinetic} code COGENT handles these difficulties by adopting a mapped multiblock discretization scheme that can enable high-accuracy and computationally efficient simulations in X-point geometries.   

\section{The Chombo Framework}
Developed at the Lawrence Berkeley National Laboratory, the DOE-supported Chombo \cite{Chombo} framework has been a \added{``}\deleted{"}developer's toolbox" for implementing highly scalable block-structured adaptive mesh refinement (AMR) algorithms for Department of Energy applications since its initial release in 2000. Because many of the algorithmic building blocks for such algorithms are both complex to implement and are also shared across many applications, Chombo provides a set of tools for implementing finite difference and finite volume methods for the solution of partial differential equations on block-structured adaptively refined logically-rectangular grids. 

Both elliptic and time-dependent modules are included. Chombo supports calculations in complex geometries with both embedded boundaries and mapped grids, and also supports particle methods. Most parallel platforms are supported, and cross-platform self-describing file formats are included via an interface with HDF5. While not formally part of xSDK \cite{xSDK}, Chombo adheres to all xSDK mandatory policies other than policy {\it M1} (the requirement for a Spack installation), relying instead on a well-documented customized build system. 

Chombo maintains interoperability with a number of community libraries, including PETSc \cite{petsc} and hypre \cite{hypre} for linear and nonlinear solvers, SUNDIALS \cite{sundials1, sundials2} for time integration, and FFTW \cite{fftw} for Fast Fourier Transforms.   Chombo has also long had an in-house regression testing system; a set of machines are dedicated to constantly check-out, build, and run a suite of Chombo tests and applications in a variety of configurations. As capabilities were added to Chombo in support of COGENT, specific regression tests were added to the Chombo suite to protect the new capabilities, which has allowed the Chombo team to quickly catch otherwise difficult-to-find changes which might impact these capabilities. 

\section{Higher-order Mapped Multiblock Finite-Volume Discretizations in Chombo}
Like many DOE problems, the physics in tokamaks have natural coordinate alignments, as described in section \ref{sect:gyrokinetic}. Taking advantage of these alignments when constructing discretizations yields substantial improvements in efficiency and accuracy, but they are often too complex for traditional mapped-grid approaches.
High-order Mapped-Multiblock (HO MMB) finite volume discretizations, on the other hand, are a general and mathematically consistent way to compose sets of localized mappings to span complex domains.

\begin{figure}[t]
\includegraphics[width=2.2in]{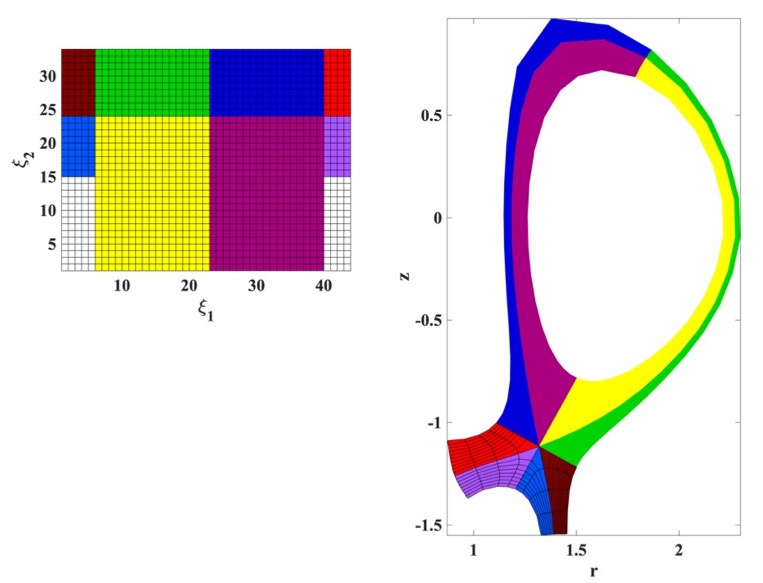}
\hspace{0.2in}
\includegraphics[width=2.2in]{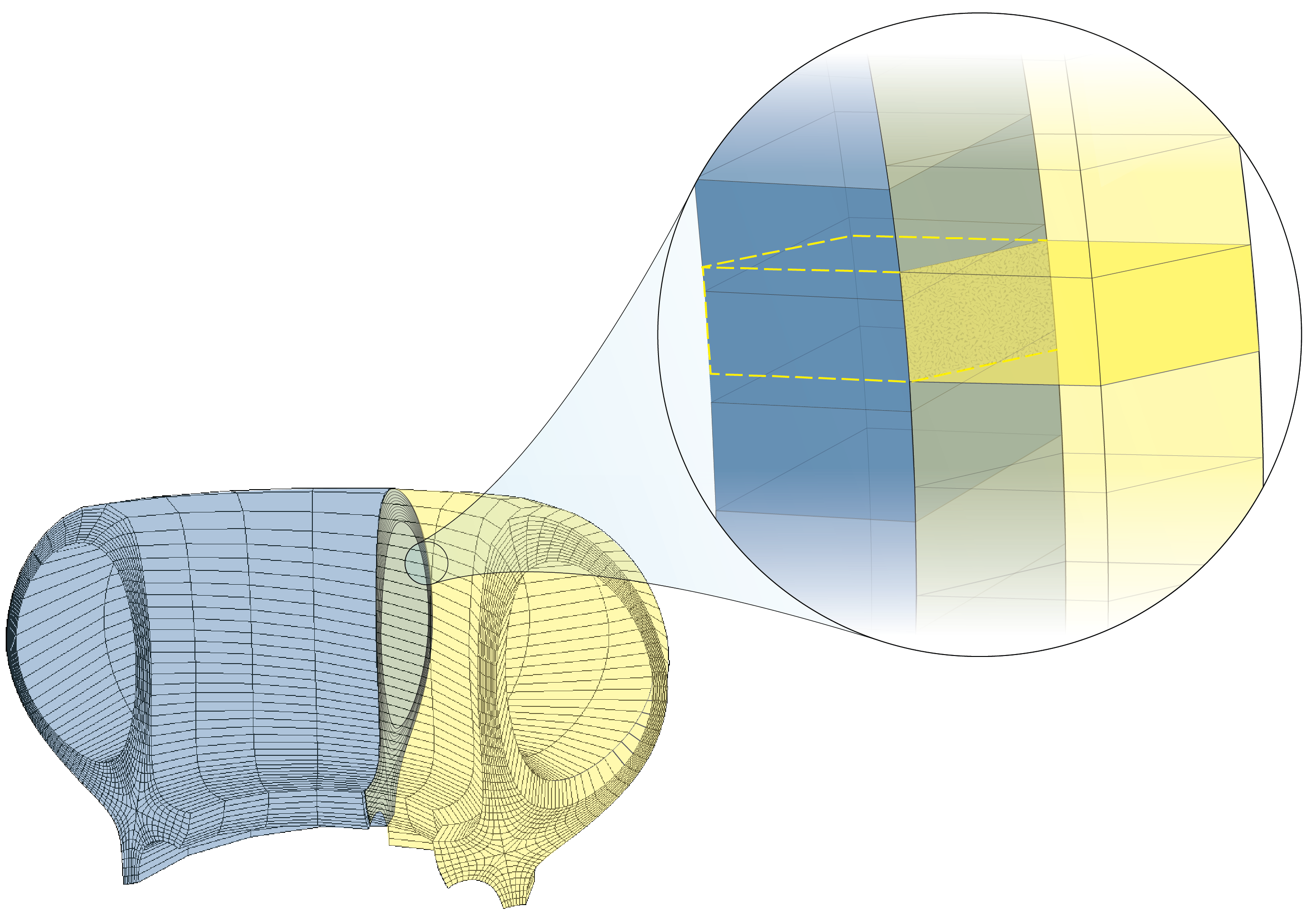}

\caption{Mapped-multiblock decomposition of the edge region of a tokamak. {\it (left)} Data is organized in logically-rectangular blocks in computational space, which corresponds to the mapped-grid cross-sectional domain (middle) (each color represents a separate mapping block, which employs a distinct mapping function. Mapping blocks meet at {\it multiblock interfaces}. (right) Full 3D mapped-multiblock discretization, showing block boundaries in the axial direction. While block boundaries are conformal in the cross-sectional mapping, we relax this requirement in the axial direction to account for the non-integral rotation of the magnetic field-aligned coordinates in the axial direction. }
\label{fig:COGENTMMB}       
\end{figure}

Our approach to these geometrically-complex computational domains is to decompose the domain into distinct regions in the configuration space in which simpler mappings can be accurately generated, and then compose these sets of \added{``}\deleted{"}mapping blocks" to represent the entire domain (see figure \ref{fig:COGENTMMB}). High-order discretizations (greater than 2) are essential for this because we expect to lose an order of accuracy at {\it multiblock boundaries}, the location of which will be fixed. (This is unlike the case for AMR, where we also lose an order of accuracy at coarse-fine interfaces, but have the flexibility to place those interfaces and their corresponding loss of accuracy in locations where the solution accuracy won't be significantly degraded. In the mapped-multiblock approach, we will not have that flexibility.) 

While it is well-known in numerical analysis that one can drop an order of accuracy on a set that is one dimension less than the domain and still retain the overall accuracy of the method,  dropping to first-order in spatial accuracy in a globally second-order method is unacceptable due to the resulting degradation in solution quality at first order. However, dropping accuracy from fourth-order to third-order is acceptable for these problems.

Our multiblock discretizations are composed of two parts:  the basic fourth-order mapped-grid discretizations in each mapping block, and the multiblock interpolation used to tie the discretizations and solutions in each block together.

\subparagraph*{Fourth-order mapped-grid discretizations}
In the interiors of mapping blocks, we base our discretizations on the fourth-order freestream preserving finite-volume formalism described in \cite{HOmapped}. Each mapping block is defined on a logically-rectangular subset of the domain, and each block is integrated using the same global timestep. The task for each of these mapping blocks is to generate mappings which are accurate and differentiable enough to compute metric terms to $O(\Delta x^5)$. For a fourth-order finite-volume method, some care must also be taken in discriminating between fourth-order point values (used for computing local physics, for example), and fourth-order cell-averaged (for computing conservative finite-volume updates) and face-averaged (for computing finite-volume fluxes) values. The resulting algorithms wind up convolving and deconvolving between averaged and point values depending on the specific needs of the algorithm.  

In software, the expression of this is object-oriented classes (derived from an interface class which defines the API), along with a set of utility functions to manage fourth-order convolution with the mapping terms and the transitions back and forth between point values and cell- and face-averaged values. 

\begin{figure}
    \centering
    \includegraphics[width=2.9in]{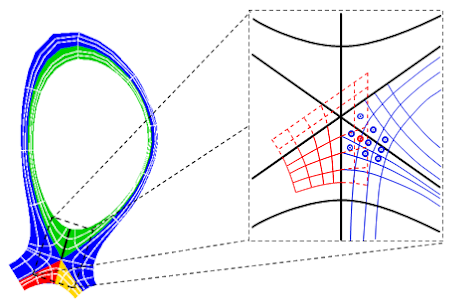}
    \caption{Example of multiblock interpolation is used as a boundary condition at multiblock boundaries.  To fill ghost cells in the red mesh, we extend its mapping beyond the MMB boundary to create ghost regions in that mapping. We then use a higher-order least-squares interpolation from cells in the blue regions to fill ghost cell values like the one shown (red circle).}
    \label{fig:mmbInterpolation}
\end{figure}

\subparagraph*{Multiblock interpolation}
To fill stencils at multiblock boundaries, we extend the mapping of each block past its boundaries to create layers of \added{``}\deleted{"}ghost cells", and then use a high-order least-squares approach to interpolate values from the relevant other mappings to the appropriately-centered locations in the ghost regions. In some cases, specific mappings may even be periodic in particular directions and wind up copying from other parts of the same mapping block. (see figure \ref{fig:mmbInterpolation}).

Developing the infrastructure required to manage the sorts of connectivity required for complex mapped-multiblock configurations required \deleted{some}\added{considerable} development work in the Chombo framework. Fortunately, Chombo's modular and extensible object-oriented design made these extensions possible. For example, we were able to use C++ inheritance to derive multiblock specializations of the classes which manage \added{interactions}\deleted{things} like connectivity between patches in the domain (which was already present to support Chombo's AMR capabilities).  This also enabled MMB capabilities to be added to the existing Chombo release versions rather than requiring deep changes and special branches.  This in turn made it easier for other applications to take advantage of these capabilities, as described briefly in Section \ref{sect:otherApps}.

\section{Mixed-dimensional Chombo}
In Chombo, dimensionality is a compile-time parameter. This is useful because many of the discretizations and operators implemented in Chombo have specific features depending on the dimensionality of the problem: for example, two-dimensional hyperbolic advection and elliptic-operator discretizations are much simpler than their three-dimensional counterparts, particularly near multiblock, domain, and AMR boundaries. Chombo also makes extensive use of multigrid approaches for its linear and nonlinear solvers -- the coarsening and refinement operations inherent in these methods are also inherently dimensional (coarsening and refining an $N\times N$ mesh in 2D is not the same as an $N \times N \times 1$ mesh in 3D, for example). The Chombo build system then encodes information such as dimensionality in the names of its libraries and executables, allowing for multiple builds of Chombo applications to co-exist in the same directories, which is often convenient for developers and users. (For instance, building an executable ``driver" built with {\tt DIM =2}, using {\tt mpicc} and mpif90 for C++ and Fortran compilers, and with optimization turned on, would produce an executable file named {\tt driver2d.Linux.64.mpiCC.mpif90.OPT.MPI.ex}, and object and dependency files generated during the build are maintained in similarly-coded subdirectories to the build directory.  This is often convenient for developers and users who wish to maintain multiple builds in a single location.)

In the phase-space discretization used by COGENT, we evolve distributions in a 4- or 5-dimensional space formed by a tensor product of a two- or three-dimensional spatial mesh (location in physical space) and a two-dimensional mesh in velocity space mesh (the combination of the gyrokinetic approximation and magnetic-field-aligned coordinates allows for the elimination of the third velocity dimension).  While there are parts of the algorithm which are carried out in the full high-dimensional space, other parts of the algorithm are carried out entirely in either physical (configuration) space or in velocity space. These parts of the algorithm are fundamentally lower-dimensional, and it is convenient and efficient to use the corresponding dimensionally-consistent version of Chombo. In particular, we want to be able to directly leverage the commonly-used two- and three-dimensional implementations in the Chombo framework where possible, rather than rewrite custom implementations which only operate on subsets of a high-dimensional space.  This enhances usability and code correctness, since the standard implementations are heavily used and tested over a range of applications, which would be less true for custom implementations used only by a single application.

To enable this, Chombo implemented a framework which supports mixed-dimensional algorithms which can use the native Chombo implementations in each dimensionality. This capability was enabled using C++ namespaces.  The ability to encapsulate Chombo in a {\tt Chombo} namespace was already required by the xSDK policies (M9) to allow for safe interoperability with other libraries \cite{xSDK}; implementing mixed-dimension programming entailed adding additional dimension-specific namespaces ({\tt Chombo::D2, Chombo::D3}, etc). A set of interface functions live outside the dimensional namespaces and manage interdimensonal operations like slicing higher dimensional spaces into lower ones and injection of data from lower-dimensional spaces into higher ones.  Another set of operators was then implemented in order to manage spreading and reduction operations in Chombo's distributed MPI domains. 

For example, imagine a phase-space algorithm with a distribution function $\phi^n$ defined over the full 4D phase space at time $t^n$: $\phi(\vec{x}, \vec{u}, t^n) = \phi(x, y, u, v, t^n)$, in which one computes a source term (like a collision operator) in velocity space (dimensions 2 and 3), integrates the source term over velocity space and applies it to a set of operations in configuration space (dimensions 0 and 1): 
\begin{enumerate}
  \item   $\phi^n  =  \phi(x,y,u,v, t^n) \ \  \rm{on} \ \  \Omega_4$
  \item   $S_u(x,y,u,v,t) = f(u,v) \ \ \rm{on} \ \ \Omega_{vel}$
  \item $S_x^n(x,y) = \int_{u,v} S(x,y,u,v,t^n) d\vec{u} $
  \item $\phi_x(x,y)^n = \int_{u,v} \phi(x,y,u,v,t^n) d\vec{u}$
  \item $\phi_x^{n+1}(x,y) = \phi_x^n(x,y) + \Delta t S_x^n(x,y)$
\end{enumerate}
Implemented in mixed-dimensional Chombo, this algorithm would look like:
\begin{enumerate}
    \item {\tt Chombo::D4}: $\phi^n  =  \phi(x,y,u,v, t^n) \ \  \rm{on} \ \  \Omega_4$
  \item   {\tt Chombo::D4}: $S_u(x,y,u,v,t) = f(u,v) \ \ \rm{on} \ \ \Omega_{vel}$
  \item {\tt Chombo}: $S_x^n(x,y) = {\tt Slice}_{4D\rightarrow2D}\bigl(\int_{u,v} S(x,y,u,v,t^n) d\vec{u} \bigr) $
  \item {\tt Chombo::D2} $\phi_x(x,y)^n = \int_{u,v} \phi(x,y,u,v,t^n) d\vec{u}$
  \item {\tt Chombo::D2} $\phi_x^{n+1}(x,y) = \phi_x^n(x,y) + \Delta t S_x^n(x,y)$
\end{enumerate}
In other words, steps 1 and 2 are computed in the full 4-dimensional space in the namespace {\tt Chombo::D4}, and steps 4 and 5 are computed in Chombo's 2-dimensional namespace {\tt Chombo::D2}. Step 3 exists in a nondimensional {Chombo} namespace, and combines a discrete reduction operator (the integral over velocity space, which is performed in {\tt Chombo::D4}) with {\tt Slice}, which is an interdimensional operator which copies a 2D slice of the 4D Chombo domain into the 2-dimensional space. 

Both the reduction and slicing operations (along with their spreading and injection operations) are designed to operate efficiently and correctly for any parallel decomposition of the domain and for the cell-centered and face-centered data centerings that are common in finite-volume algorithms. Implementing the mixed-dimensional Chombo build required significant extensions to the Chombo build system to incorporate the correct system of namespaces and library and application source-code compilation. In practice, however, it greatly simplifies implementing phase-space algorithms such as that used by COGENT because it allows developers to leverage Chombo implementations designed specifically for the dimensionality that one is operating in.

\section{COGENT}
The finite-volume gyrokinetic code COGENT (COntinuum Gyrokinetic Edge New Technology) has been developed by the Edge Simulations Laboratory (ESL) collaboration for edge plasma modeling. The code has a number of physical models including (i) a fully-kinetic approach in which a 5D gyrokinetic equation is employed for all plasma species (i.e., electrons and ions) and coupled to a 3D electrostatic gyro-Poisson equation for electrostatic potential perturbations $\Phi(R)$ \cite{Lee2018}, and (ii) a hybrid approach in which a gyrokinetic equation for the ion species is coupled to a 3D quasi-neutrality equation for the vorticity variable, $\varpi=\nabla_\perp(\alpha \nabla_\perp \Phi)$, and a 3D fluid model for the electron species \cite{Dorf2021}. The gyrokinetic equation represents an advection equation, where the phase-space velocity depends on the values of electric field, ${\bf E}=-\nabla \Phi$. When collisions are included, drag and diffusion terms in the velocity space are added. The gyro-Poisson equation is an elliptic equation for a perpendicular (to the magnetic field) Laplacian operator, and the quasi-neutrality vorticity equation is written for a time derivative of a vorticity variable and involves both stiff (diffusion-like) and non-stiff terms in the right-hand-side \cite{Dorf2021}. 

The ion-scale transport simulations (which is the main focus of the COGENT code) are stiff and could benefit greatly from implicit time integration in order to step over fast time scales associated, for example, with electron dynamics. In the case of a hybrid model, all fast time scales are contained within the 3D field/fluid part of the hybrid system and therefore only the 3D system needs to be treated implicitly. The 5D gyrokinetic ion system can be treated explicitly. To that end, a consistent high-order time integration approach that includes implicit treatment of selected stiff terms is implemented in COGENT. It is based on semi-implicit additive Runge–Kutta (ARK) methods and employs the Jacobian-free Newton-Krylov (JFNK) approach to handle nonlinearities. Multigrid solvers from hypre are used for preconditioning purposes to improve JFNK convergence properties.

The COGENT code has axisymmetric (2D+2V) and non-axisymmetric (3D+2V) versions, with the former version used to assess the properties of collisional transport and the latter version employed to study the effects of microturbulence. Either collisional or microturbulence transport is highly anisotropic and requires the use of mapped multiblock technology to perform simulations in X-point geometries.

Development of COGENT benefited a great deal from being able to take advantage of Chombo's high-order mapped-multiblock and mixed-dimension programming capabilities (which were developed in anticipation of COGENT's needs). COGENT developers were able to focus on the specific algorithmic needs of building an edge-plasma gyrokinetic code without having to worry about building specific infrastructure for these capabilities. In many cases, COGENT developers were able to use existing Chombo capabilities in combination with HOMMB or mixed-dimension support rather than write custom implementations.  Close coordination between COGENT developers and the Chombo development team ensured that issues which did arise were able to be addressed through cooperative efforts.

\subsection{Example results}

\begin{figure}
    \centering
    \includegraphics[width=2.75in]{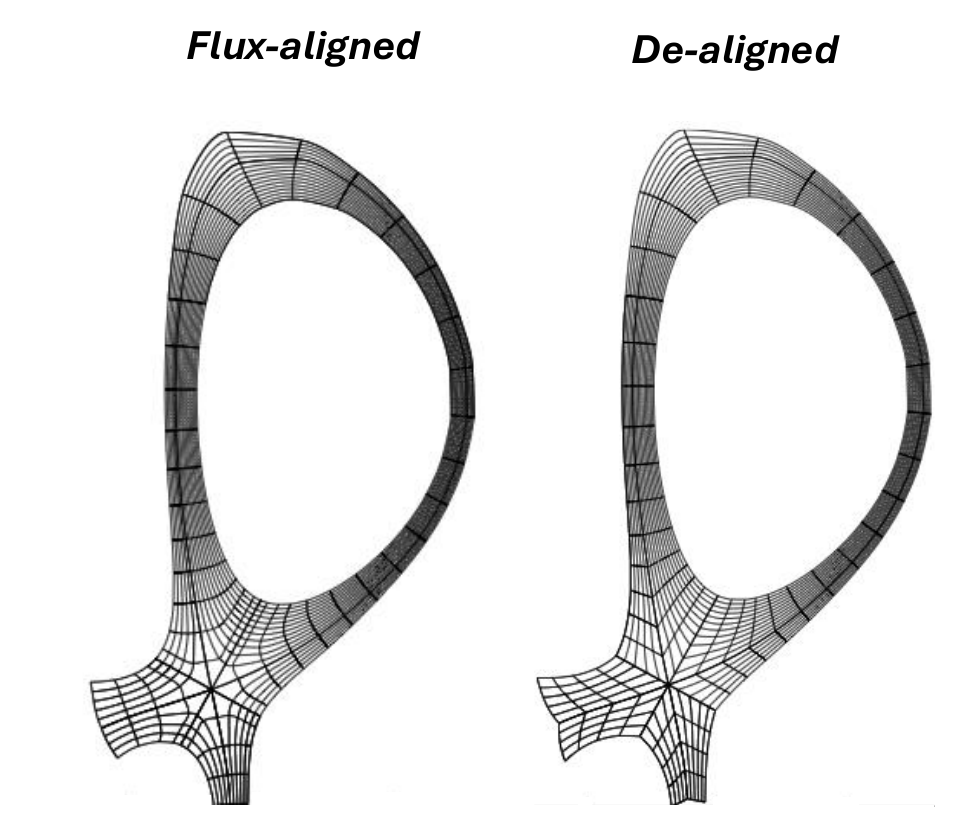}
    \caption{Generation of a multiblock mapping grid. The mesh on the right has poloidal grid lines that coincide with level surfaces of a modified flux function. A modified flux function is obtained blending the original flux function (that sourced the fully flux-aligned grid on the left) with a block-aligned linear function near the X-point. }
    \label{fig:alignedMeshes}
\end{figure}  

A successful application of high-order mapped multiblock discretization was demonstrated for axisymmetric ion gyrokinetic advection in a single-null (X-point) geometry \cite{Dorr2018}. The poloidal grid lines of a global grid are obtained as level surfaces of a modified flux function, $\Psi$, given by blending the original flux function, $\psi_0$ (takes the form of $\Psi_0 - \Psi_X \approx \alpha^2 R^2 - \beta^2 Z^2 $ in the vicinity of the X point) and the block-aligned linear function $\psi_{lin} = D (|\alpha R| - |\beta Z|)$. As a result, the grid is rectilinear and block-aligned in the vicinity of the X-point and aligned with the original flux function outside the transition radial distance $D$ (see Figure \ref{fig:alignedMeshes}). A test case of a Boltzmann equilibrium demonstrated fourth-order convergence in the block interiors, and third-order convergence near the block boundaries where interpolation is used \cite{Dorr2018}. 
\begin{figure}
    \centering
    \includegraphics[width=2.75in]{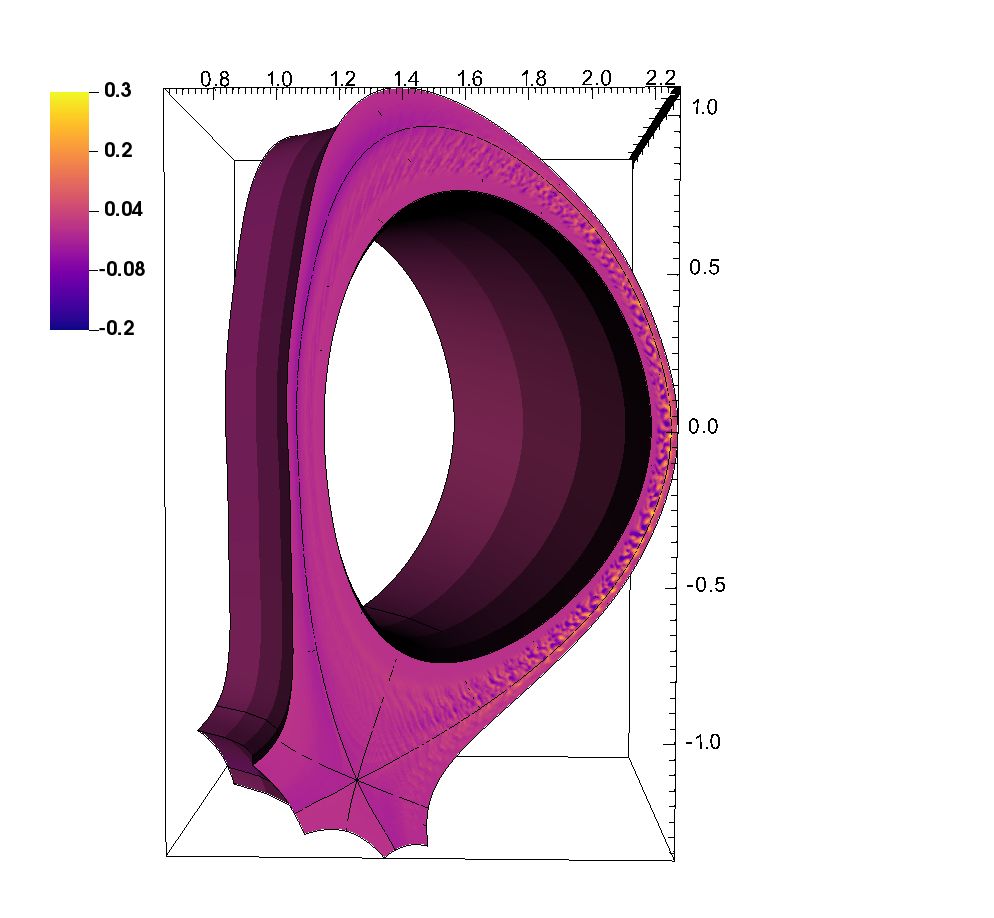}
    \caption{COGENT simulations of ion-scale drift microturbulence obtained with the hybrid 5D gyrokinetic ion - 3D fluid electron model. A magnetic geometry of the DIII-D tokamak and model plasma parameters are used for this example simulation, which includes $3x10^8$ phase-space cells. Shown are perturbations of electrostatic potential. Black curves illustrate poloidal block boundaries.}
    \label{fig:movieFrame}
\end{figure}

Given that transport anisotropy is mitigated near the X-point, abandoning the original magnetic flux surface alignment in that region can be consistent with tolerable numerical pollution for ion advection. Accordingly, successful physics studies of axisymmetric ion transport were performed \cite{Dorf2016}. However, electron transport has a much higher degree of anisotropy and is more sensitive to numerical pollution. As a result, it is less amenable to grid de-alignment. COGENT and Chombo teams are currently working to overcome this challenge by combining the MMB approach with mesh refinement near the X-point to handle electron transport. For the case of a hybrid simulation model, mesh refinement only needs to be applied to the low-dimensional fluid electron system (in the configuration space), whereas the high-dimensional kinetic ion system could be evolved on the original de-aligned grid without additional refinement. 

While this development is still ongoing, the COGENT physics team were able to perform the first-ever continuum axisymmetric \cite{Dorf2016, Dorf2018} and non-axisymmetric \cite{Dorf2021, Dorf2022} edge transport simulations in a single-null geometry by making use of the MMB approach combined with fully flux-aligned meshes. While the use of such grids minimizes numerical pollution, it comes at the expense of degraded convergence properties \added{because these problems involve diverging metric factors at an X-point}. An example of a ion scale microturbuelnce simulation in a realistic DIII-D geometry is shown in Figure \ref{fig:movieFrame}. To exploit the strong anisotropy of microturbulence, the 3D control volumes are both flux-aligned and field-aligned (as described in Figure \ref{fig:COGENTMMB}).

\section{Other applications \label{sect:otherApps}}
One of the major advantages of software frameworks like Chombo is the ability to leverage development efforts across multiple \deleted{applicaotions}\added{applications}. While the Chombo capabilities developed in this work were specifically designed to support COGENT, they have found use in other applications as well. As an example, we will describe their use in a space weather application.

\subsection{Space Weather}
\begin{figure}
    \centering
    \includegraphics[width=2.75in]{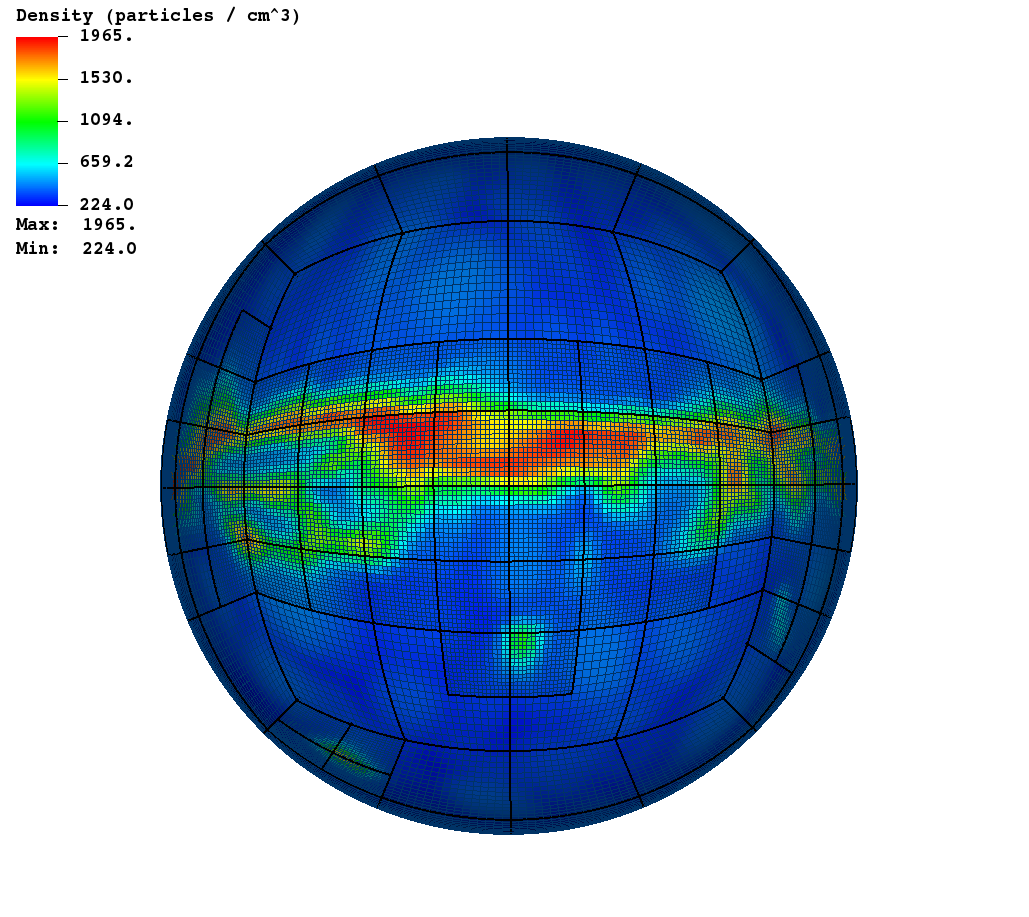}
    \caption{HelioCubed simulation showing the plasma density at the inner boundary of the Cubed Sphere domain, which is placed at 21.5 solar radii. The plasma density at this distance is described by the semi-empirical Wang--Sheeley--Arge corona model. Thin blue and thicker black lines represent the cell and patch boundaries, respectively. AMR is used to obtain finer resolution near the equatorial current sheet. Two levels of refinement are shown. Figure courtesy of Christopher Bozhart and Talwinder Singh.}
    \label{fig:spaceWeather}
\end{figure}

One example where the higher-order mapped multiblock support developed for COGENT is proving useful is the HelioCubed space weather modeling effort led by Nikolai Pogorolev at the University of Alabama, Huntsville (UAH) in collaboration with Chombo developers at LBNL\cite{Singh2023, Pogorelov2024} The UAH effort has used the Chombo framework for many years for related space weather modeling efforts.  The natural coordinates to use for these calculations is spherical coordinates with the sun as the center.  However, computing in spherical coordinate systems suffers from the presence of coordinate singularities at the polar axis. In particular, the global time steps in explicit schemes are often driven by the stability requirement determined by exceedingly small computational cells near the poles. By switching to a cubed-sphere coordinate system using the Chombo mapped-multiblock support initially developed for COGENT, the HelioCubed effort has demonstrated an ability to increase their stable timesteps by two orders of magnitude (from 6 seconds to 600 seconds), while also maintaining high spatial accuracy due to the $4^{th}$-order discretizations employed in Chombo's MMB infrastructure.  An example calculation from this effort is shown in figure \ref{fig:spaceWeather}.

\added{HelioCubed was able to take advantage of Chombo's object-oriented implementation of mapped-multiblock geometry support.  Implementing a new mapping requires implementation of a minimal set of geometry-specific functions in a derived class, which is located with the application code.  The functions which convolve and deconvolve face-averaged, cell-averaged, and point values are implemented as a part of the Chombo MMB support, and interact with the geometry specifications provided by the application in a modular way.
The space-weather application uses a cubed-sphere mapping, a common mapping used in a wide range of applications. Because of this, we had implemented a cubed-sphere mapping along with several example mappings of various complexity  for testing purposes as we developed the MMB infrastructure.  As a result, the space-weather geometry required no real extra effort to implement. However, specifying a new geometry unique to an application requires only the subset of functions in the geometry derived class which specify the mapping itself.}

\section{Looking to the Future -- Proto, AMR, embedded boundaries, and sparse grids}
While the Chombo/COGENT partnership embodies a successful example of how frameworks can support applications, progress and improvements are always underway. The Chombo framework developers continue to support the COGENT effort through a number of capability improvements, which we expect will continue to prove useful to other applications as well.  Examples of these improvements include:
\begin{enumerate}
    \item \textbf{Proto for GPU performance portability}
    Performance portability and GPU support for
the higher-order mapped-multiblock efforts is provided via Proto \cite{Proto}, the ECP-supported Chombo performance-portability layer and DSL.  Support for the HOMMB infrastructure has been extended to Proto, and key computational kernels in  COGENT will be ported to GPUs using Proto to demonstrate the resulting performance improvements.
    \item \textbf{Adaptive Mesh Refinement for HOMMB}
    Current COGENT simulation campaigns have made clear specific needs for localized higher spatial resolution around the X-point, demonstrating the importance of adaptive mesh refinement (AMR) as a critical tool for HOMMB. This is particularly true in higher-dimensional computations due their extreme use of computational resources for even moderate resolutions. However, ensuring correct discretizations which preserve higher-order accuracy at coarse-fine interfaces in the presence of MMB meshes requires some care. Prototype implementations of AMR for HOMMB for Poisson’s equation are currently being extended to support the more complex equation sets required by COGENT and will soon be coupled to COGENT for refinement around the X-point in tokamaks.
    \item \textbf{Higher-order embedded-boundary discretizations for HOMMB}
    Finally, even with AMR and HOMMB meshes, there are often still complex geometrical features not
aligned with coordinates (e.g. divertor plates in tokamaks, planets in space weather) which must be accurately represented. We are currently extending the higher-order embedded-boundary cut-cell approaches developed by the Chombo team \cite{HOEB} to represent these accurately in HOMMB discretizations.
    \item \textbf{Sparse Grids} While coordinate alignment and AMR somewhat mitigate the curse of dimensionality that inherently plagues the five-dimensional phase-space simulations COGENT performs, these simulations remain highly resource intensive.  The COGENT and Chombo teams are currently exploring the use of sparse grid methods \cite{SGreview} to further mitigate the curse of dimensionality.  This necessitates both modifications to Chombo's native stencils to retain 4$^\textrm{th}$-order accuracy in the presence of sparse grids' dependence on high-order mixed derivatives and a more flexible MPI-communicator infrastructure to leverage the additional parallelism enabled by sparse grids.  We expect these modifications will also be useful for other Chombo applications, particularly coupled multiphysics applications and those looking to explore parallel-in-time integration schemes.
\end{enumerate}

\section{Conclusion}
As described above, the specific needs of the COGENT tokamak edge-modeling effort required a set of new capabilities from the Chombo framework in order to support a set of novel equations and discretizations.  Due to Chombo's flexible object-oriented design, the Chombo development team was able to work closely with the COGENT development team to design and implement a fairly complex set of new features, including high-order mapped-multiblock discretizations and mixed-dimensional programming.  These features were essential to the COGENT team's success, and in turn fed back into the main Chombo development and release pipeline, allowing other development efforts to take advantage of these capabilities. This level of support also required close collaboration from the mathematicians, software developers, and physicists involved in this effort. 

\section*{ \added{Code Availability} } 
\added{COGENT may be obtained from the COGENT GitHub repository: \\ {\tt https://github.com/LLNL/COGENT}. 
\\
Chombo documentation and instructions may be found at {\tt https://Chombo.lbl.gov} and the Chombo Github repository: \\ {\tt https://github.com/applied-numerical-algorithms-group-lbnl/Chombo\_3.2}  }

\section*{Acknowledgments}
Work at Berkeley Lab was supported by the Director, Office of Science, of the U.S. Department of Energy and the Advanced Scientific Computing Research program under Contract No. DE-AC02-05CH11231. Work at Lawrence Livermore National Laboratory was supported by the U.S. Department of Energy under contract DE-AC52-07NA27344.

\end{document}